\documentclass[12pt, a4paper]{amsart}
\usepackage{amsmath,amssymb,amsthm,booktabs,longtable}
\usepackage{enumitem}
\usepackage{comment}
\usepackage{graphicx}
\usepackage[all]{xy}
\usepackage{caption}
\usepackage{color}
\usepackage[colorlinks=true, linkcolor=blue, citecolor=blue, urlcolor=blue]{hyperref}
\hypersetup{breaklinks=true}
\usepackage{mathdots}
\usepackage{url}
\usepackage[T1]{fontenc} 
\usepackage{here}
\usepackage{lscape}
\usepackage{tikz}

\theoremstyle{plain}
\makeatletter 

\def\l@section{\@tocline{1}{0pt}{1em}{}{}}
\def\l@subsection{\@tocline{2}{0pt}{2em}{5em}{}}
\def\l@subsubsection{\@tocline{3}{0pt}{3em}{5em}{}}

\@addtoreset{table}{section}

\@addtoreset{figure}{section}
\makeatother %
\newtheorem{theorem}{Theorem}[section]

\newtheorem{corollary}[theorem]{Corollary}
\newtheorem{proposition}[theorem]{Proposition}

\newtheorem{problem}[theorem]{Problem}
\newtheorem{question}[theorem]{Question}

\theoremstyle{definition}
\newtheorem{definition}[theorem]{Definition}

\newtheorem{remark}[theorem]{Remark}

\newcommand{\Aut}{\mathop{\mathrm{Aut}}\nolimits}

\newcommand{\R}{\mathbb{R}}
\newcommand{\C}{\mathbb{C}}
\newcommand{\N}{\mathbb{N}}

\newcommand{\HA}{\mathbb{H}}

\makeatletter
    
    \@addtoreset{equation}{section}
  \makeatother
\makeatother

\title[]{Deformations of
Standard Locally Homogeneous Spaces}
\author{Kazuki Kannaka and Toshiyuki Kobayashi}


\address[Kazuki KANNAKA]{%
Faculty of Mathematics and Physics, Institute of Science and Engineering, Kanazawa University, Kakumamachi, Kanazawa, Ishikawa, 920-1192, JAPAN;
RIKEN Interdisciplinary Theoretical and Mathematical Sciences (\lowercase{i}THEMS), 
Wako, Saitama 351-0198, Japan.
}
\email{kannaka@se.kanazawa-u.ac.jp}

\address[Toshiyuki KOBAYASHI]{%
Graduate School of Mathematical Sciences,
The University of Tokyo, 3-8-1 Komaba, Meguro, Tokyo, 153-8914, Japan;
French-Japanese Laboratory in Mathematics and its Interactions,
FJ-LMI CNRS IRL2025, Tokyo, Japan
}
\email{toshi@ms.u-tokyo.ac.jp}

\date{\today}

\begin{document}

\begin{abstract}
Let $X=G/H$ be a homogeneous space,
where $G \supset H$ are reductive Lie groups.
We ask: in the setting where $\Gamma \backslash G/H$ is a standard quotient, to what extent can the discrete subgroup $\Gamma$ be deformed while preserving
the proper discontinuity of the $\Gamma$-action on $X$?

We provide several classification results, including: conditions under which local rigidity holds for compact standard quotients $\Gamma\backslash X$; criteria for when a standard quotient can be deformed into a nonstandard one; 
a characterization of the maximal Zariski-closure of discontinuous groups under small deformations; and conditions under which Zariski-dense deformations occur.

Proofs of the results stated in this paper are provided in detail in 
arXiv:2507.03476.
\end{abstract}

\maketitle%


\section{Standard Locally Homogeneous Spaces}

We consider a homogeneous space $X=G/H$, where
 $G \supset H$ are  
real algebraic groups. 
A \emph{discontinuous group} for $X$ is a discrete subgroup $\Gamma \subset G$ that acts properly discontinuously and freely on $X$.
Any $G$-invariant local geometric structure on $X$ then descends to the quotient space $X_\Gamma:=\Gamma \backslash X$, which is Hausdorff and constitutes a typical example of a $(G,X)$-manifold in the sense of Ehreshmann and Thurston.
We investigate the following question:
\begin{problem}
\label{mainquestion} 
To what extent can cocompact discontinuous groups $\Gamma$ for $X$ be deformed? 
\end{problem}

In the classical setting where $X$ is an irreducible Riemannian symmetric space,
the Selberg--Weil local rigidity theorem (\cite{Weil_discrete_subgroups})
asserts that cocompact discontinuous groups $\Gamma$ admit no nontrivial deformations,  with the sole exception of the two-dimensional case, where their deformations (up to equivalence) are parametrized by 
 the Teichm\"uller space.
 
In contrast to classical rigidity, a notable phenomenon arising when $H$ is noncompact---first observed by the second author in the early 1990s \cite{Kobayashi93}---is that cocompact discontinuous groups for \emph{pseudo-Riemannian} symmetric spaces $X=G/H$ exhibit
greater \lq\lq flexibility\rq\rq. 
Specifically, there exist arbitrarily high-dimensional compact quotients $X_\Gamma$ of irreducible symmetric spaces $X=G/H$ that admit continuous deformation (\cite[Thms.~A and B]{Kobayashi98}).
Nevertheless, 
it has also been shown that 
certain compact quotients $X_\Gamma$ in pseudo-Riemannian settings exhibit local rigidity (\cite[Prop.\ 1.8]{Kobayashi98}).

The difficulty regarding Problem~\ref{mainquestion} 
 is that when 
 $H$ is not compact, small deformations of a discrete subgroup can easily destroy the proper discontinuity.
 For example, in the setting of nonreductive groups, consider the cocompact discontunuous group $\Gamma \simeq \mathbb{Z}$ on $X \simeq \mathbb{R}$ by translations. A small deformation of this action within the affine transformation group generally fails to remain properly discontinuous.

Goldman~\cite{Goldmannonstandard} conjectured, in the context of the 3-dimensional compact anti-de Sitter space, that any small deformation of a \emph{standard} cocompact discontinuous group preserves proper discontinuity.
This stability conjecture was proved affirmatively in \cite{Kobayashi98},
using the properness criterion established in
\cite{Benoist96, Kobayashi89, Kobayashi96}.
The \emph{stability} of proper discontinuity was further studied in \cite{Kassel12} for certain \emph{standard quotients}.

We now recall from \cite[Def.~1.4]{KasselKobayashi16}
the notion of standard quotients $X_\Gamma$:

\begin{definition}
\label{def:standard}
A discontinuous group $\Gamma$ for $X$ is called \emph{standard} if there exists a reductive subgroup $L \subset G$ which contains $\Gamma$ and acts properly on $X$. 
\end{definition}

When $L \subset G$ acts properly and cocompactly on $X=G/H$, any torsion-free, cocompact, discrete subgroup $\Gamma \subset L$ (which exists by a theorem of Borel) yields a compact, locally homogeneous space $X_\Gamma$.

 A criterion for a triple $(G,H,L)$ of reductive Lie groups to ensure that $L$ acts properly on $X=G/H$ was established in \cite[Thm.~4.1]{Kobayashi89}, along with a cocompactness criterion in \cite[Thm.~4.7]{Kobayashi89}. 
 These criteria are computable and yield a list of such triples $(G,H,L)$, which we present in Table~\ref{tab:cpt-CK-sym}.
Here, $L_{ss}$ denotes the semisimple part of $L$.
Entries without a dash correspond to cases in which $G/H$ is a symmetric space; these are taken from  \cite{KobayashiYoshino05}. 
Entries marked with a dash indicate cases obtained by interchanging the roles of $H$ and $L$.
Table~\ref{tab:cpt-CK-sym}
 completes the earlier list presented in
Kulkarni~\cite{Kulkarni1981proper} and the second author~\cite{Kobayashi89, Kobayashi_kawaguchiko90,Kobayashi1997discontinuouspseudoRiemannian}. 
In particular, Cases~4 and 5 are from
 \cite{Kulkarni1981proper}, and Cases~1, 1', 2, 3, 4', and 5' are from \cite{Kobayashi89}.
Recent work by Tojo~\cite{Tojo2019Classification} and
by Boche\'{n}ski--Tralle~\cite{BochenskiTralle24} contributes to the supporting evidence for the (essential) completeness of Table~\ref{tab:cpt-CK-sym}.

\section{Deformation of Discontinuous Groups}

Let $\Gamma$ be a finitely generated discrete subgroup of a Lie group $G$.
By a deformation, we mean fixing the abstract group structure of $\Gamma$, and varying the homomorphisms from $\Gamma$ into the group $G$.

We now return to the main theme of this article: the deformation of \emph{discontinuous groups} for homogeneous spaces $X=G/H$, particularly
when $H$ is noncompact.
In such cases, a discrete subgroup of $G$ may not act properly discontinuously on $G/H$, and
small deformations can easily destroy proper discontinuity. 
To capture this phenomenon, the second author introduced in \cite{Kobayashi93}  
 the \emph{deformation space} $\mathcal{R}(\Gamma,G;X)$, defined as the set of faithful homomorphisms $\varphi\colon \Gamma \to G$ such that $\varphi(\Gamma)$ acts properly discontinuously on $X$.

There are two natural actions on  $\mathcal{R}(\Gamma,G;X)$: one by the automorphism group
$\Aut(\Gamma)$, and the other by the inner automorphism group of $G$.
These actions commute with each other.

\begin{definition}
\label{def:locally_rigid_for_X}
(\cite[Sect.\ 1]{Kobayashi98}).
We say that $\varphi \in \mathcal{R}(\Gamma, G; X)$ is \emph{locally rigid as a discontinuous group for $X$}, or that the quotient space $X_{\varphi(\Gamma)}$ is \emph{locally rigid}, 
if the $G$-orbit through $\varphi$ is open in $\mathcal{R}(\Gamma, G; X)$ with respect to the compact-open topology.
\end{definition}

When $H$ is compact, this notion coincides with the classical notion of local rigidity in Weil~\cite{Weil_remark}.

To better understand Definition~\ref{def:locally_rigid_for_X},
we recall from \cite[Sect.~5.3]{Kobayashi-unlimit} the notions of 
the \emph{higher Teichm\"uller space}
$\mathcal{T}(\Gamma,G;X)$ 
and the \emph{moduli space}
$\mathcal{M}(\Gamma,G;X)$, defined
as follows:
\begin{align*}
\mathcal{T}(\Gamma,G;X) &:= \mathcal{R}(\Gamma,G;X)/G.
\\
\mathcal{M}(\Gamma,G;X) &:= \Aut(\Gamma)\backslash \mathcal{R}(\Gamma,G;X)/G.
\end{align*}
These spaces capture what may be regarded as \emph{essential deformations} of discontinuous groups when $G$ acts faithfully on $X$.
Note that the term \emph{higher Teichm\"uller space} is typically used in the setting where $\Gamma$ is a surface group and $X = G$, whereas here we consider a broader context in which $X = G/H$ is a homogeneous space with $H$ possibly noncompact.

\begin{table}[htb!]
\begin{tabular}{c|c|c}
Case & $X=G/H$ & $L_{ss}$ \\
\hline
1 & $SU(2n,2)/Sp(n,1)$ & $SU(2n,1)$
\\
\hline
1'-1 & $SU(2n,2)/SU(2n,1)$ $(n\geq 2)$ & $Sp(n,1)$ \\
1'-2 &
$SU(2,2)/SU(2,1)$ & $Spin(4,1)$ \\
\hline
2-1 & $SU(2n,2)/U(2n,1)$ $(n\geq 2)$ & $Sp(n,1)$ \\
2-2 &
$SU(2,2)/U(2,1)$ & $Spin(4,1)$ \\
\hline
3 & 
$SO(2n,2)/U(n,1)$ & $SO(2n,1)$\\
\hline
4-1 &
$SO(2n,2)/SO(2n,1)$ $(n\geq 2)$ & 
$SU(n,1)$ \\
4-2 &
$SO(2,2)/SO(2,1)$ & 
$SU(1,1)$ \\
\hline
4' & 
$SO(2n,2)/SU(n,1)$ & $SO(2n,1)$ \\
\hline
5-1 &
$SO(4n,4)/SO(4n,3)$ $(n\geq 2)$ & 
$Sp(n,1)$ \\
5-2 &
$SO(4,4)/SO(4,3)$ & 
$Spin(4,1)$ \\
\hline
5' &
$SO(4n,4)/Sp(n,1)$ & 
$SO(4n,3)$ \\
\hline
6 &
$SO(8,8)/SO(8,7)$ & 
$Spin(8,1)$ \\
\hline
6' &
$SO(8,8)/Spin(8,1)$ & 
$SO(8,7)$ \\
\hline
7 &
$SO(4,4)/(SO(4,1)\times SO(3))$ &
$Spin(4,3)$ \\
\hline
7' &
$SO(4,4)/Spin(4,3)$ &
$SO(4,1)$ \\
\hline
8 &
$SO(8,\mathbb{C})/SO(7,\mathbb{C})$ &
$Spin(7,1)$ \\
\hline
8' &
$SO(8,\mathbb{C})/Spin(7,1)$ &
$SO(7,\mathbb{C})$  \\
\hline
9 &
$SO(8,\mathbb{C})/SO(7,1)$ &
$Spin(7,\mathbb{C})$ \\
\hline
9' &
$SO(8,\mathbb{C})/Spin(7,\mathbb{C})$ &
$SO(7,1)$  \\
\hline
10 &
$SO^*(8)/U(3,1)$ &
$Spin(6,1)$  \\
\hline
10' &
$SO^*(8)/Spin(6,1)$ &
$SU(3,1)$  \\
\hline
11 &
$SO^*(8)/
(SO^*(6)\times SO^*(2)) $ & $Spin(6,1)$ \\
\hline
11' &
$SO^*(8)/Spin(6,1)$ & $SO^*(6)$\\
\hline
12 &
$SO(4,3)/(SO(4,1)\times SO(2))$ 
& $G_{2(2)}$  \\
\hline
12' &
$SO(4,3)/G_{2(2)}$ 
& $SO(4,1)$  
\end{tabular}
\caption{$L$ acts properly and cocompactly on $G/H$}
\label{tab:cpt-CK-sym}
\end{table}

\section{Classification Results}
With this background, we are ready to provide a rigorous formulation of Problem~\ref{mainquestion} concerning deformations of \emph{standard cocompact discontinuous groups}, and to present classification results based on Table~\ref{tab:cpt-CK-sym}.

    \begin{question}
\label{question:deform-ck}
Classify all triples $(G, H, L)$ listed in Table~\ref{tab:cpt-CK-sym}
that admit a torsion-free, cocompact discrete subgroup $\Gamma \subset L$ satisfying each of the following properties: 
\begin{enumerate}[label=$(Q\arabic*).$]
    \item 
    \label{item:local-rigid}
    $\Gamma$ is not locally rigid as a discontinuous group for $X=G/H$     (see Definition~\ref{def:locally_rigid_for_X});
    \item 
    \label{item:non-standard}
    $\Gamma$ can be deformed into a nonstandard discontinuous group for $X$ (see Definition~\ref{def:standard});
    \item
    \label{item:Zariski-dense}
    $\Gamma$ can be deformed into a Zariski-dense subgroup of $G$, while preserving proper discontinuity of the action on $X$.
\end{enumerate}
\end{question}

Regarding $(Q3)$, we note that when $H$ is noncompact, the properness condition forces $L\neq G$, which in turn implies that a standard discontinuous group $\Gamma$ cannot be Zariski-dense in $G$. 

We now state the classification results.
\begin{theorem}
\label{thm:classification}
   Let $(G,H,L)$ be one of the triples listed in Table~\ref{tab:cpt-CK-sym}.
 \newline
\noindent (1).
Condition $(Q1)$ holds
    if and only if the case is one of the following: 1, 1'-2, 2-2, 3, 4, 4', 5-2, 7', 10, 10', 11, 12'. 
 \newline\noindent(2). 
 Conditions $(Q2)$ and $(Q3)$ are equivalent, and they hold if and only if the case is one of the following: 1'-2, 2-2, 3, 4-2, 4', 5-2, 7', 10, 11, 12'.
\end{theorem}

Moreover, the result is invariant under (appropriately defined) local isomorphisms of the triple $(G,H,L)$, which is also part of the theorem.
\begin{remark}
 The implication $(Q2) \Rightarrow (Q3)$ may fail when $G$ is semisimple rather than simple.
\end{remark}
\begin{remark}
Prior results include:
\newline\noindent {(1)}. 
Nontrivial deformations 
are known to exist 
in Cases 1, 4-1, 4-2, 7', 12' of Table~\ref{tab:cpt-CK-sym} \cite{Kobayashi98};
\newline\noindent {(2)}. 
Local rigidity 
in Cases 2-1, 5-1, 5' \cite{Kobayashi98};
\newline\noindent{(3)}. 
Nonstandard deformations for group manifolds
$(G\times G)/\Delta(G)$ with $G=SO(n,1)$ or $SU(n,1)$ \cite{Goldmannonstandard, Kobayashi98};
\newline\noindent {(4)}. 
Zariski-dense deformations 
in Case 3 \cite{Kassel12}. 
\end{remark}
\begin{remark}
Spectral analysis on pseudo-Riemannian locally symmetric spaces is a developing research area.  Discrete spectrum under deformations of $X_\Gamma$ is studied in \cite{KasselKobayashi16}, while spectral decomposition for the standard $X_\Gamma$ is developed in \cite{KasselKobayashi2019standard}.
\end{remark}

Theorem~\ref{thm:classification} yields notable new results, e.g.: 
\begin{corollary}
    There exists a compact, seven- dimensional space form $\Gamma\backslash SO(4,4)/SO(3,4)$
    of signature $(4,3)$ and negative curvature,
    such that $\Gamma$ is Zariski-dense in the full isometry group $SO(4,4)$.
    \end{corollary}
\begin{corollary} 
Every compact standard locally symmetric space
\newline $\Gamma\backslash SO(8,\C)/SO(7,\C)$
and $\Gamma\backslash SO(8,8)/SO(7,8)$ is locally rigid.
\end{corollary}

We outline the proof of Theorem~\ref{thm:classification}. Upper bounds on deformations are discussed in Section~\ref{sec:upper_bound}, and the construction of deformations in Section~\ref{sec:bending}.
\section{Upper Bounds of Deformations}
\label{sec:upper_bound}
 
The upper bound on deformations---including, as a special case, local rigidity of discontinuous groups (Definition~\ref{def:locally_rigid_for_X})---is obtained via infinitesimal methods. 

For $(Q1)$, a sufficient condition for local rigidity as a discontinuous group can be derived, following \cite{Kobayashi98}, from the local rigidity of a homomorphism $\iota\colon \Gamma \hookrightarrow L \rightarrow G$, based on Weil's vanishing condition for $H^1(\Gamma, \mathfrak{g})$ \cite{Weil_remark}.
This vanishing can be verified using Raghunathan's theorem \cite{Raghunathan65}.

For $(Q2)$, to show that nontrivial deformations exist but no nonstandard ones occur---that is, $(Q1)$ holds but $(Q2)$ fails---we prove the following proposition assuming vanishing of cohomology with smaller coefficients in a setting where $H^1(\Gamma, \mathfrak{g}) \neq 0$. 
\begin{proposition}
    \label{prop:local-rigid-g/l}
 Let $G$ be a Lie group, and $L$ its closed subgroup.
  If $\Gamma \subset L$ is a finitely presented subgroup satisfying $H^1(\Gamma, \mathfrak{g}/\mathfrak{l}) = 0$, then  $\Gamma$
  admits no nontrivial deformations outside $L$ up to $G$-conjugacy.
\end{proposition}
Proposition~\ref{prop:local-rigid-g/l} applies to
Case~4-1 ($n\geq 3$) and Case~10'.
However, Theorem~\ref{thm:classification} also treats subtler cases (such as Case~1 and Case~4-1 for $n=2$), where $L \simeq U(n,1)$ and 
 $H^1(\Gamma, \mathfrak{g}/\mathfrak{l}) \neq 0$
at the infinitesimal level, yet this does not lift to actual local deformations.
In these cases, the proof utilizes the idea of 
Goldman--Millson~\cite{goldman-millson87}
(see also \cite{Klingler11}), which involves
the cup product in cohomology.

\section{Bending Construction for \texorpdfstring{$Spin(n,1)$}{Spin(n,1)} }
\label{sec:bending}
We now outline the method to construct deformations.  
Previously known approaches include one using the nonvanishing of $H^1(\Gamma, \mathbb{R})$, as in \cite{Kobayashi98}, for $\mathfrak{l}=\mathfrak{so}(n,1)$ or $\mathfrak{su}(n,1)$, and another
using the bending construction by Johnson--Millson~\cite{JoMi84}, as applied in Kassel~\cite{Kassel12} for $L = SO(n,1)$.

Some homogeneous spaces $X$, such as the 15-dimensional space form $SO(8,8)/SO(7,8)$ (see Case 6 of Table~\ref{tab:cpt-CK-sym}),
admit a proper action by reductive subgroups only when the group is globally isomorphic to $Spin(n,1)$. In particular, $SO(n,1)$ cannot act properly on such spaces.
Moreover, 
not every cocompact discrete subgroup of $SO(n,1)$ can be lifted to $Spin(n,1)$ for $n \ge 4$
(Martelli--Riolo--Slavich~\cite{nonspin-hyperbolic}).

To complete the proof of Theorem~\ref{thm:classification}, we extend the bending construction to cocompact discrete subgroups of $Spin(n,1)$.
This formulation allows us to cover cases where $\Gamma \subset SO(n,1)$ as well.
We also determine the upper bound of the Zariski closure of the deformed image (Theorem~\ref{thm:bending}).

Reformulating Millson's theorem \cite{Millson76} in terms of Clifford algebras, we can construct arithmetic subgroups $\Gamma$ with specific geometric properties:
\begin{proposition}
    \label{proposition:Millson-spin}
    For every $k\in \N$, 
    there exists a cocompact discrete subgroup $\Gamma \subset Spin(n,1)$ such that the associated compact hyperbolic manifold 
    $M:=\Gamma\backslash \HA^{n}$ contains
   $k$ mutually disjoint, orientable, connected, totally geodesic closed hypersurfaces $N_{1},\ldots, N_{k}$
with the property that the complement $S:=M\smallsetminus (N_{1}\cup\cdots \cup N_{k})$ is connected.
\end{proposition}

    \begin{figure}
        \centering
        \includegraphics[width=5cm]{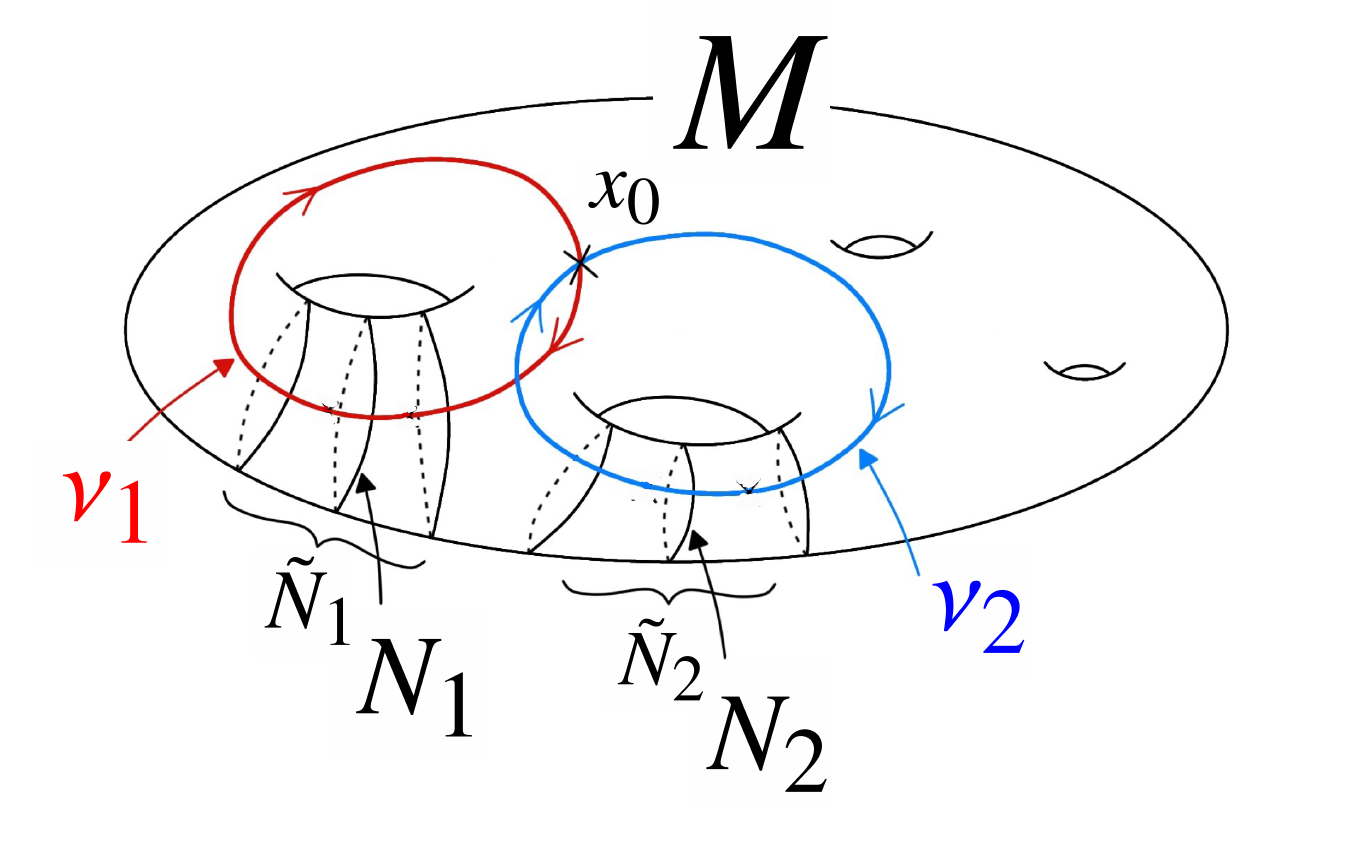}, 
        \caption{Loops $\nu_{1}$ and $\nu_{2}$ in the case $k=2$.}
        \label{fig:hnn2}
    \end{figure}

A deformation of
this arithmetic subgroup $\Gamma$ can be constructed inside a Lie group $G$ via a homomorphism $\varphi\colon Spin(n,1) \to G$,
using the bending construction developed by Johnson--Millson~\cite{JoMi84}.
To this end, we choose oriented loops $\nu_i$ for $i=1,\ldots,k$, each intersecting $N_i$ transversally exactly once and  disjoint from the others (see Fig.~\ref{fig:hnn2}).
We realize $\pi_{1}(N_i)$  as subgroups of $G$ via the composition: 
\[
\pi_{1}(N_{i})\hookrightarrow
\pi_{1}(M)=\Gamma\subset Spin(n,1)
\xrightarrow{\varphi}G,
\]
and make use of the fact that the group structure of $\Gamma$ can be described as an iterated HNN extension of $\pi_1(S)$ to derive the following proposition.

\begin{proposition}
\label{prop:Johnson-Millson}
Choose $v_{i}\in\mathfrak{g}$ to be fixed under the action of $\operatorname{Ad}\circ \varphi (\pi_{1}(N_{i}))$
for each $i=1,\ldots,k$.
Then there exists a unique one-parameter family of homomorphisms
$\varphi_{t}\colon \Gamma\rightarrow G$
such that
$\varphi_{t}\equiv
\varphi$ on $\pi_{1}(S)$ and 
$\varphi_{t}([\nu_{i}])=\varphi([\nu_{i}])\exp(tv_{i})$
for $1 \le i \le k$. 
\end{proposition}

The geometric condition on
$v_{i}\in \mathfrak{g}$ required in Proposition~\ref{prop:Johnson-Millson} is ensured by a representation-theoretic condition: namely, that $v_i$ is fixed by the subgroup associated with the
totally geodesic hypersurface $N_{i}$.
This subgroup appears in $Spin(n,1)$ as one isomorphic to $Spin(n-1,1)$.

The \emph{spherical harmonics representations} of degree $j$ are real irreducible finite-dimensional representations $V_j$ $(j=0,1,\ldots)$ of $\mathfrak{spin}(n,1)$ with a nonzero $\mathfrak{spin}(n-1,1)$-fixed vector.
Each $V_j$ arises as a subrepresentation of $C^\infty(S^n)$ consisting of eigenfunctions of  $\Delta_{S^{n}}$ with eigenvalues $-j(j+n-1)$.

Let $G$ be a Zariski-connected real algebraic group, and $\varphi\colon Spin(n,1) \to G$ a nontrivial homomorphism, where $n \ge 2$.
Let $\mathfrak{g}^{\varphi}$ denote the Lie subalgebra of $\mathfrak{g}$, generated by 
$d\varphi(\mathfrak{spin}(n,1))$ together with all $d \varphi(\mathfrak{spin}(n,1))$-submodules that are isomorphic to irreducible
spherical harmonics representations of $\mathfrak{spin}(n,1)$.

Let $G^\varphi$ denote the Zariski-connected algebraic subgroup 
 with Lie algebra $\mathfrak{g}^{\varphi}$ (such exists).
\begin{theorem}
    \label{thm:bending}
{(1)}
        \label{introitem:realization-gphi}
        There exists a cocompact arithmetic subgroup $\Gamma \subset Spin(n,1)$ such that $\varphi(\Gamma)$
    can be deformed so that its Zariski-closure is  
    $G^{\varphi}$.
    \label{introitem:upper-bound}
\newline\noindent{(2)}
    The group $G^{\varphi}$ gives the upper bound for the Zariski-closure of $\varphi'(\Gamma)$, up to $G$-conjugacy, for any small deformation $\varphi'$ of $\varphi|_\Gamma$, provided that $n \geq 3$.
\end{theorem}

Theorem~\ref{thm:bending} guarantees the existence of Zariski-dense deformations, provided that
$\mathfrak{g}^{\varphi}=\mathfrak{g}$.

By contrast, Theorem~\ref{thm:bending} also applies to
the special case in which $\mathfrak{g}^\varphi$ is minimal, leading to another notable corollary.
Let $[\mathfrak{g}:V_j]$ denote the multiplicity of
the $\mathfrak{spin}(n,1)$-module $V_j$ occurring in $\mathfrak{g}$ under the adjoint action via $\varphi$.
This multiplicity vanishes for all sufficiently large $j$.

    \begin{corollary}
    \label{cor:deform_to_centralizer}
    In the setting of Theorem~\ref{thm:bending} with $n \ge 3$,
 the following conditions are equivalent:
\begin{enumerate}[label=(\roman*)]
\item
$[\mathfrak{g}:V_j]=0$ for all $j \ge 1$. 
\item
Every torsion-free cocompact discrete subgroup $\Gamma \subset L=Spin(n,1)$ admits no nontrivial deformations beyond $\varphi(L)\cdot Z_{G}(\varphi(L))$.
\end{enumerate}
\end{corollary}
\begin{remark}
Klingler~\cite[Thm.~1.3.7]{Klingler11} established a sufficient condition for an analogous condition to (ii) of Corollary~\ref{cor:deform_to_centralizer}
when $L=SU(n,1)$, .
\end{remark}

Theorem~\ref{thm:bending}~(2) is proved using Proposition~\ref{prop:local-rigid-g/l} and Raghunathan's vanishing theorem (\cite{Raghunathan65}).

For the proof of (1) of Theorem~\ref{thm:bending},
we take sufficiently large $k$
and choose appropriate vectors $v_{1},\ldots, v_{k}\in \mathfrak{g}$ in 
Propositions~\ref{proposition:Millson-spin} and~\ref{prop:Johnson-Millson}.

An optimal value of $k$ is given by
\begin{equation}
    \label{def:k}
    k:=\max(2, [\mathfrak{g}:V_{1}], 
    [\mathfrak{g}:V_{2}],\ldots)  \in \N.
    \end{equation}
The optimality of \eqref{def:k} is shown using the general results presented in the next section.

\section{Family of Zariski-Dense Subgroups}
We introduce three invariants $\underline{\eta}(G) \le \eta(G) \le \overline\eta(G)$ as the minimal cardinalities of finite subsets $\mathcal{X} \subset \mathfrak{g}$ satisfying the following conditions, respectively:
\newline\
     $\overline\eta(G): \
      G(t\mathcal{X})=G  \textrm{ for all } t>0$;
\newline\
$        \eta(G): \
        G(t\mathcal{X})=G
        \textrm{ for any sufficiently small }t>0$; 
\newline\
  $      \underline{\eta}(G): \
       G(\mathcal{X})=G$.
\newline
Here, $G(\mathcal{X})$ denotes the Zariski closure in $G$ of the subgroup generated by $e^{X}$ for all $X\in \mathcal{X}$. 

This definition makes sense without reductiveness assumption of $G$.  For instance, for
    $G=\R^{n}$ regarded as a unipotent  commutative algebraic group,
    we have
 $\underline{\eta}(G)= \eta(G)=\overline{\eta}(G)=n$.

\begin{remark} For semisimple $G$,
it is known that $\underline{\eta}(G) = 2$; see, e.g., 
\cite[Lem.~5.3.13]{Labourie_lecture} for a proof when $G$ is the split group $SL(n, \mathbb{R})$. 
However, the invariants $\eta$ and $\overline\eta$ are designed to control the behavior of the family $G(t \mathcal{X})$ as the parameter $t \in \R$ varies.
\end{remark}

\begin{theorem}
\label{thm:two_generates_reductive}
    Let $G$ be a Zariski-connected real reductive algebraic group. Then we have
    \[
    \underline{\eta}(G)= \eta(G)=\overline\eta(G)=2,
    \]
    except for the following cases: 
    \begin{itemize}
        \item ($\underline{\eta}, \eta, \overline{\eta})=(1,1,1)$
        for a split $\R$-torus $G$;
\item $(\underline\eta, \eta, \overline\eta)=(1,2,2)$ for 
        a nonsplit $\R$-torus $G$; 
\item $(\underline\eta, \eta, \overline\eta)=(2,2,3)$ when 
$\mathfrak{su}(2)$ is 
an ideal of
        $\mathfrak{g}$.
    \end{itemize}
\end{theorem}

\medskip
Detailed proofs of the results presented here, as well as certain deformation results in the noncompact case of $X_\Gamma$ extending \cite{KannakaOkudaTojo24},
will appear in \cite{KannakaKobayashi25}.

\section*{Acknowledgments}
The first author was partially supported by Special Postdoctoral
Researcher Program at RIKEN and JSPS 
Kakenhi 
(JP24K16929).
The second author was partially supported by JSPS Kakenhi
(JP23H00084)
 and by the Institut Henri Poincar\'e (Paris) and the Institut des Hautes {\'E}tudes Scientifiques.
 (Bures-sur-Yvette).  

\end{document}